\documentclass{article}
\usepackage{amsfonts}
\usepackage{amsmath}
\usepackage{amssymb}

\input{tcilatex}

\begin{document}

\title{$L_{p}$ Convergence with Rates of Smooth Poisson-Cauchy Type Singular
Operators}
\author{George A. Anastassiou \& Razvan A. Mezei}
\date{$\quad$}
\maketitle



\begin{centering}
Department of Mathematical Sciences\\
The University of Memphis\\
Memphis, TN 38152, U.S.A.\\
{\tt ganastss@memphis.edu}\\
{\tt rmezei@memphis.edu}\\
\end{centering}

\pagestyle{myheadings} 
\markboth{Rates of Smooth Gauss-Weierstrass Singular
Operators}{George A. Anastassiou, Razvan A. Mezei}

\begin{quote}
\textbf{Abstract.} In this article we continue the study of smooth
Poisson-Cauchy Type singular integral operators on the line regarding their
convergence to the unit operator with rates in the $L_{p}$ norm, $p\geq 1$.
The related established inequalities involve the higher order $L_{p}$
modulus of smoothness of the engaged function or its higher order derivative.


\end{quote}

\section{Introduction}

The rate of convergence of singular integrals has been studied in [9], [13],
[14], [15], [7], [8], [4], [5], [6] and these articles motivate our work.
Here we study the $L_{p}$, $p\geq 1$, convergence of smooth Poisson-Cauchy
Type singular integral operators over $\mathbb{R}$ to the unit operator with
rates over smooth functions with higher order derivatives in $L_{p}(\mathbb{R%
}).$ We establish related Jackson type inequalities involving the higher $%
L_{p}$ modulus of smoothness of the engaged function or its higher order
derivative. The discussed operators are not in general positive, see [10],
[11]. Other motivation comes from [1], [2].

\section{Results}

In the next we introduce and deal with the \textit{smooth Poisson-Cauchy
Type singular integral operators} $M_{r,\xi }(f;x)$ defined as follows.

For $r\in \mathbb{N}$ and $n\in \mathbb{Z}_+$ we set

\noindent 
\begin{equation}
\alpha _{j}=\left\{ 
\begin{array}{ll}
(-1)^{r-j}{\binom{r}{j}}j^{-n}, & \quad j=1,\ldots ,r, \\ 
1-\dsum\limits_{j=1}^{r}(-1)^{r-j}{\binom{r}{j}}j^{-n}, & \quad j=0,%
\end{array}%
\right.
\end{equation}%
that is $\sum\limits_{j=0}^{r}\alpha _{j}=1$.

Let $f\in C^{n}(\mathbb{R})$ and $f^{(n)}\in L_{p}(\mathbb{R}),$ $1\leq
p<\infty $, $\alpha \in \mathbb{N},$ $\beta >\frac{1}{2\alpha },$ we define
for $x\in \mathbb{R}$, $\xi >0$ the Lebesgue integral 
\begin{equation}
M_{r,\xi }(f;x)=W\int_{-\infty }^{\infty }\frac{\sum_{j=0}^{r}\alpha
_{j}f(x+jt)}{\left( t^{2\alpha }+\xi ^{2\alpha }\right) ^{\beta }}dt,
\end{equation}%
where the constant is defined as%
\begin{equation*}
W=\frac{\Gamma \left( \beta \right) \alpha \xi ^{2\alpha \beta -1}}{\Gamma
\left( \frac{1}{2\alpha }\right) \Gamma \left( \beta -\frac{1}{2\alpha }%
\right) }.
\end{equation*}

$\bigskip $

\textbf{Note 1. }The operators $M_{r,\xi }$ are not, in general, positive.
See [10], (18).

We notice by $W\int_{-\infty }^{\infty }\frac{1}{\left( t^{2\alpha }+\xi
^{2\alpha }\right) ^{\beta }}dt=1$, that $M_{r,\xi }(c,x)=c$, $c$ constant,
see also [10], [11], and 
\begin{equation}
M_{r,\xi }(f;x)-f(x)=W\left( \sum_{j=0}^{r}\alpha _{j}\int_{-\infty
}^{\infty }\left[ f(x+jt)-f(x)\right] \frac{1}{\left( t^{2\alpha }+\xi
^{2\alpha }\right) ^{\beta }}dt\right) .
\end{equation}%
We use also that 
\begin{equation}
\int_{-\infty }^{\infty }\frac{t^{k}}{\left( t^{2\alpha }+\xi ^{2\alpha
}\right) ^{\beta }}dt=\left\{ 
\begin{array}{ll}
0, & \quad k\ \hbox{odd}, \\ 
\frac{1}{\xi ^{2\alpha \beta -k-1}\alpha }\frac{\Gamma \left( \frac{k+1}{%
2\alpha }\right) \Gamma \left( \beta -\frac{k+1}{2\alpha }\right) }{\Gamma
\left( \beta \right) }, & \quad k\ \hbox{even,}\text{ with }\beta >\frac{k+1%
}{2\alpha },%
\end{array}%
\right.
\end{equation}%
see [16].

We need the $r$th $L_{p}$-modulus of smoothness 
\begin{equation}
\omega _{r}(f^{(n)},h)_{p}:=\sup_{\left\vert t\right\vert \leq h}\Vert
\Delta _{t}^{r}f^{(n)}(x)\Vert _{p,x},\ \ h>0,
\end{equation}%
where 
\begin{equation}
\Delta _{t}^{r}f^{(n)}(x):=\sum_{j=0}^{r}(-1)^{r-j}{\binom{r}{j}}%
f^{(n)}(x+jt),
\end{equation}%
see [12], p.~44. Here we have that $\omega _{r}(f^{(n)},h)_{p}<\infty $, $%
h>0 $.

We need to introduce 
\begin{equation}
\delta _{k}:=\sum_{j=1}^{r}\alpha _{j}j^{k},\quad k=1,\ldots ,n\in \mathbb{N}%
,
\end{equation}%
and denote by $\lfloor \cdot \rfloor $ the integral part. Call 
\begin{equation}
\tau (w,x):=\sum_{j=0}^{r}\alpha _{j}j^{n}f^{(n)}(x+jw)-\delta
_{n}f^{(n)}(x).
\end{equation}%
Notice also that 
\begin{equation}
-\sum_{j=1}^{r}(-1)^{r-j}{\binom{r}{j}}=(-1)^{r}{\binom{r}{0}}.
\end{equation}%
According to [3], p.\ 306, [1], we get 
\begin{equation}
\tau (w,x)=\Delta _{w}^{r}f^{(n)}(x).
\end{equation}%
Thus 
\begin{equation}
\Vert \tau (w,x)\Vert _{p,x}\leq \omega _{r}(f^{(n)},|w|)_{p},\quad w\in 
\mathbb{R}.
\end{equation}

Using Taylor's formula, and the appropriate change of variables, one has
(see [6])%
\begin{equation}
\sum_{j=0}^{r}\alpha _{j}[f(x+jt)-f(x)]=\sum_{k=1}^{n}\frac{f^{(k)}(x)}{k!}%
\delta _{k}t^{k}+\mathcal{R}_{n}(0,t,x),
\end{equation}%
where 
\begin{equation}
\mathcal{R}_{n}(0,t,x):=\int_{0}^{t}\frac{(t-w)^{n-1}}{(n-1)!}\tau
(w,x)dw,\quad n\in \mathbb{N}.
\end{equation}%
Using the above terminology we obtain for $\beta >\frac{2\left\lfloor \frac{n%
}{2}\right\rfloor +1}{2\alpha }$ that 
\begin{equation}
\Delta (x):=M_{r,\xi }(f;x)-f(x)-\sum_{m=1}^{\lfloor n/2\rfloor }\frac{%
f^{(2m)}(x)\delta _{2m}}{(2m)!}\frac{\Gamma \left( \frac{2m+1}{2\alpha }%
\right) \Gamma \left( \beta -\frac{2m+1}{2\alpha }\right) }{\Gamma \left( 
\frac{1}{2\alpha }\right) \Gamma \left( \beta -\frac{1}{2\alpha }\right) }%
\xi ^{2m}=\mathcal{R}_{n}^{\ast }(x),
\end{equation}%
where 
\begin{equation}
\mathcal{R}_{n}^{\ast }(x):=W\int_{-\infty }^{\infty }\mathcal{R}_{n}(0,t,x)%
\frac{1}{\left( t^{2\alpha }+\xi ^{2\alpha }\right) ^{\beta }}dt,\quad n\in 
\mathbb{N}.
\end{equation}%
In $\Delta (x)$, see (14), the sum collapses when $n=1$.

We present our first result.

\medskip \noindent \textbf{Theorem 1.} \textit{Let} $p,q>1$ \textit{such
that }$\frac{1}{p}+\frac{1}{q}=1$, $n\in \mathbb{N},$ $\alpha \in \mathbb{N}%
, $ $\beta >\frac{1}{\alpha }\left( \frac{1}{p}+n+r\right) $\textit{\ and
the rest as above. Then} 
\begin{equation}
\Vert \Delta (x)\Vert _{p}\leq \frac{\left( 2\alpha \right) ^{\frac{1}{p}%
}\Gamma \left( \beta \right) \Gamma \left( \frac{q\beta }{2}-\frac{1}{%
2\alpha }\right) ^{\frac{1}{q}}\xi ^{n}\tau ^{\frac{1}{p}}}{\Gamma \left( 
\frac{q\beta }{2}\right) ^{\frac{1}{q}}\Gamma \left( \frac{1}{2\alpha }%
\right) ^{\frac{1}{p}}\Gamma \left( \beta -\frac{1}{2\alpha }\right) \left(
rp+1\right) ^{\frac{1}{p}}\left[ (n-1)!\right] (q(n-1)+1)^{1/q}}\omega
_{r}(f^{(n)},\xi )_{p},
\end{equation}%
\textit{where} 
\begin{equation}
0<\tau :=\left[ \int_{0}^{\infty }\left( 1+u\right) ^{rp+1}\,\frac{u^{np-1}}{%
\left( u^{2\alpha }+1\right) ^{p\beta /2}}du-\int_{0}^{\infty }\,\frac{%
u^{np-1}}{\left( u^{2\alpha }+1\right) ^{p\beta /2}}du\right] <\infty .
\end{equation}%
\textit{Hence as }$\xi \rightarrow 0$ \textit{we obtain }$\Vert \Delta
(x)\Vert _{p}\rightarrow 0$.

If additionally $f^{(2m)}\in L_{p}(\mathbb{R}),m=1,2,\ldots ,\left\lfloor 
\frac{n}{2}\right\rfloor $ then $\left\Vert M_{r,\xi }(f)-f\right\Vert
_{p}\rightarrow 0,$ as $\xi \rightarrow 0.$

\smallskip \noindent \textbf{Proof.} We observe that 
\begin{eqnarray}
|\Delta (x)|^{p} &=&W^{p}\left\vert \int_{-\infty }^{\infty }\mathcal{R}%
_{n}(0,t,x)\frac{1}{\left( t^{2\alpha }+\xi ^{2\alpha }\right) ^{\beta }}%
dt\right\vert ^{p}  \notag \\
&\leq &W^{p}\left( \int_{-\infty }^{\infty }\left\vert \mathcal{R}%
_{n}(0,t,x)\right\vert \frac{1}{\left( t^{2\alpha }+\xi ^{2\alpha }\right)
^{\beta }}dt\right) ^{p}  \notag \\
&\leq &W^{p}\left( \int_{-\infty }^{\infty }\left\vert \int_{0}^{|t|}\frac{%
(|t|-w)^{n-1}}{(n-1)!}|\tau (sign(t)\cdot w,x)|dw\right\vert \frac{1}{\left(
t^{2\alpha }+\xi ^{2\alpha }\right) ^{\beta }}dt\right) ^{p}.
\end{eqnarray}%
Hence we have 
\begin{equation}
I:=\int_{-\infty }^{\infty }|\Delta (x)|^{p}dx\leq W^{p}\left( \int_{-\infty
}^{\infty }\left( \int_{-\infty }^{\infty }\gamma (t,x)\frac{1}{\left(
t^{2\alpha }+\xi ^{2\alpha }\right) ^{\beta }}dt\right) ^{p}dx\right) ,
\end{equation}%
where 
\begin{equation}
\gamma (t,x):=\int_{0}^{|t|}\frac{(|t|-w)^{n-1}}{(n-1)!}|\tau (sign(t)\cdot
w,x)|dw\geq 0.
\end{equation}%
Therefore by using H\"{o}lder's inequality suitably we obtain%
\begin{eqnarray}
R.H.S.(19) &=&W^{p}\left( \int_{-\infty }^{\infty }\left( \int_{-\infty
}^{\infty }\gamma (t,x)\frac{1}{\left( t^{2\alpha }+\xi ^{2\alpha }\right)
^{\beta }}dt\right) ^{p}dx\right)  \notag \\
&=&W^{p}\cdot \left( \int_{-\infty }^{\infty }\left( \int_{-\infty }^{\infty
}\gamma (t,x)\frac{1}{\left( t^{2\alpha }+\xi ^{2\alpha }\right) ^{\beta /2}}%
\frac{1}{\left( t^{2\alpha }+\xi ^{2\alpha }\right) ^{\beta /2}}dt\right)
^{p}dx\right)  \notag \\
&\leq &W^{p}\cdot \left( \int_{-\infty }^{\infty }\left( \int_{-\infty
}^{\infty }\left[ \gamma (t,x)\frac{1}{\left( t^{2\alpha }+\xi ^{2\alpha
}\right) ^{\beta /2}}\right] ^{p}dt\right) \left( \int_{-\infty }^{\infty }%
\left[ \frac{1}{\left( t^{2\alpha }+\xi ^{2\alpha }\right) ^{\beta /2}}%
\right] ^{q}dt\right) ^{\frac{p}{q}}dx\right)  \notag
\end{eqnarray}%
\begin{eqnarray}
&=&W^{p}\cdot \left( \int_{-\infty }^{\infty }\left( \int_{-\infty }^{\infty
}\gamma ^{p}(t,x)\frac{1}{\left( t^{2\alpha }+\xi ^{2\alpha }\right)
^{p\beta /2}}dt\right) \left( \int_{-\infty }^{\infty }\frac{1}{\left(
t^{2\alpha }+\xi ^{2\alpha }\right) ^{q\beta /2}}dt\right) ^{\frac{p}{q}%
}dx\right)  \notag \\
&=&W^{p}\cdot \left( \int_{-\infty }^{\infty }\left( \int_{-\infty }^{\infty
}\gamma ^{p}(t,x)\frac{1}{\left( t^{2\alpha }+\xi ^{2\alpha }\right)
^{p\beta /2}}dt\right) dx\right) \left( \frac{\Gamma \left( \frac{1}{2\alpha 
}\right) \Gamma \left( \frac{q\beta }{2}-\frac{1}{2\alpha }\right) }{\Gamma
\left( \frac{q\beta }{2}\right) \alpha \xi ^{q\alpha \beta -1}}\right) ^{%
\frac{p}{q}}  \notag \\
&=&\frac{\xi ^{p\alpha \beta -1}\alpha \left[ \Gamma \left( \beta \right) %
\right] ^{p}\Gamma \left( \frac{q\beta }{2}-\frac{1}{2\alpha }\right) ^{%
\frac{p}{q}}}{\Gamma \left( \frac{q\beta }{2}\right) ^{\frac{p}{q}}\Gamma
\left( \frac{1}{2\alpha }\right) \Gamma \left( \beta -\frac{1}{2\alpha }%
\right) ^{p}}\left( \int_{-\infty }^{\infty }\left( \int_{-\infty }^{\infty
}\gamma ^{p}(t,x)\frac{1}{\left( t^{2\alpha }+\xi ^{2\alpha }\right)
^{p\beta /2}}dt\right) dx\right) .
\end{eqnarray}%
Again by H\"{o}lder's inequality we have 
\begin{equation}
\gamma ^{p}(t,x)\leq \frac{\left( \int_{0}^{|t|}|\tau (sign(t)\cdot
w,x)|^{p}dw\right) }{((n-1)!)^{p}}\,\frac{|t|^{np-1}}{(q(n-1)+1)^{p/q}}\,.
\end{equation}%
Consequently we have 
\begin{eqnarray*}
R.H.S.(21) &\leq &\frac{\xi ^{p\alpha \beta -1}\alpha \left[ \Gamma \left(
\beta \right) \right] ^{p}\Gamma \left( \frac{q\beta }{2}-\frac{1}{2\alpha }%
\right) ^{\frac{p}{q}}}{\Gamma \left( \frac{q\beta }{2}\right) ^{\frac{p}{q}%
}\Gamma \left( \frac{1}{2\alpha }\right) \Gamma \left( \beta -\frac{1}{%
2\alpha }\right) ^{p}} \\
&&\cdot \left( \int_{-\infty }^{\infty }\left( \int_{-\infty }^{\infty }%
\frac{\left( \int_{0}^{|t|}|\tau (sign(t)\cdot w,x)|^{p}dw\right) }{%
((n-1)!)^{p}}\,\frac{|t|^{np-1}}{(q(n-1)+1)^{p/q}}\frac{1}{\left( t^{2\alpha
}+\xi ^{2\alpha }\right) ^{p\beta /2}}dt\right) dx\right) \\
&=&:(\ast ),
\end{eqnarray*}%
(calling 
\begin{equation}
c_{1}:=\frac{\xi ^{p\alpha \beta -1}\alpha \left[ \Gamma \left( \beta
\right) \right] ^{p}\Gamma \left( \frac{q\beta }{2}-\frac{1}{2\alpha }%
\right) ^{\frac{p}{q}}}{\Gamma \left( \frac{q\beta }{2}\right) ^{\frac{p}{q}%
}\Gamma \left( \frac{1}{2\alpha }\right) \Gamma \left( \beta -\frac{1}{%
2\alpha }\right) ^{p}((n-1)!)^{p}(q(n-1)+1)^{p/q}})
\end{equation}%
and 
\begin{eqnarray}
(\ast ) &=&c_{1}\left( \int_{-\infty }^{\infty }\left( \int_{-\infty
}^{\infty }\left( \int_{0}^{|t|}|\tau (sign(t)\cdot w,x)|^{p}dw\right)
\,|t|^{np-1}\frac{1}{\left( t^{2\alpha }+\xi ^{2\alpha }\right) ^{p\beta /2}}%
dx\right) dt\right)  \notag \\
&=&c_{1}\left( \int_{-\infty }^{\infty }\left( \int_{-\infty }^{\infty
}\left( \int_{0}^{|t|}|\Delta _{sign(t)\cdot w}^{r}f^{(n)}(x))|^{p}dw\right)
\,|t|^{np-1}\frac{1}{\left( t^{2\alpha }+\xi ^{2\alpha }\right) ^{p\beta /2}}%
dx\right) dt\right)  \notag \\
&=&c_{1}\left( \int_{-\infty }^{\infty }\left( \int_{-\infty }^{\infty
}\left( \int_{0}^{|t|}|\Delta _{sign(t)\cdot w}^{r}f^{(n)}(x))|^{p}dw\right)
dx\right) \,|t|^{np-1}\frac{1}{\left( t^{2\alpha }+\xi ^{2\alpha }\right)
^{p\beta /2}}dt\right)  \notag \\
&=&c_{1}\left( \int_{-\infty }^{\infty }\left( \int_{0}^{|t|}\left(
\int_{-\infty }^{\infty }|\Delta _{sign(t)\cdot
w}^{r}f^{(n)}(x))|^{p}dx\right) dw\right) \,|t|^{np-1}\frac{1}{\left(
t^{2\alpha }+\xi ^{2\alpha }\right) ^{p\beta /2}}dt\right)  \notag \\
&\leq &c_{1}\left( \int_{-\infty }^{\infty }\left( \int_{0}^{|t|}\omega
_{r}(f^{(n)},w)_{p}^{p}dw\right) \,|t|^{np-1}\frac{1}{\left( t^{2\alpha
}+\xi ^{2\alpha }\right) ^{p\beta /2}}dt\right) .
\end{eqnarray}

So far we have proved 
\begin{equation}
I\leq c_{1}\left( \int_{-\infty }^{\infty }\left( \int_{0}^{|t|}\omega
_{r}(f^{(n)},w)_{p}^{p}dw\right) \,|t|^{np-1}\frac{1}{\left( t^{2\alpha
}+\xi ^{2\alpha }\right) ^{p\beta /2}}dt\right) .
\end{equation}%
By [12], p.\ 45 we have 
\begin{equation}
(R.H.S.(25))\leq c_{1}\left( \omega _{r}(f^{(n)},\xi )_{p}\right) ^{p}\left(
\int_{-\infty }^{\infty }\left( \int_{0}^{|t|}\left( 1+\frac{w}{\xi }\right)
^{rp}dw\right) \,|t|^{np-1}\frac{1}{\left( t^{2\alpha }+\xi ^{2\alpha
}\right) ^{p\beta /2}}dt\right) =:(\ast \ast ).
\end{equation}%
But we see that 
\begin{equation}
(\ast \ast )=\left( \frac{\xi c_{1}}{rp+1}\right) \left( \omega
_{r}(f^{(n)},\xi )_{p}\right) ^{p}\mathcal{J},
\end{equation}%
where 
\begin{eqnarray}
\mathcal{J} &=&\int_{-\infty }^{\infty }\left( \left( 1+\frac{\left\vert
t\right\vert }{\xi }\right) ^{rp+1}-1\right) \,|t|^{np-1}\frac{1}{\left(
t^{2\alpha }+\xi ^{2\alpha }\right) ^{p\beta /2}}dt  \notag \\
&=&2\int_{0}^{\infty }\left( \left( 1+\frac{t}{\xi }\right) ^{rp+1}-1\right)
\,t^{np-1}\frac{1}{\left( t^{2\alpha }+\xi ^{2\alpha }\right) ^{p\beta /2}}%
dt.
\end{eqnarray}%
Here we find 
\begin{eqnarray}
\mathcal{J} &=&2\xi ^{p\left( n-\alpha \beta \right) }\int_{0}^{\infty
}\left( \left( 1+u\right) ^{rp+1}-1\right) \,u^{np-1}\frac{1}{\left(
u^{2\alpha }+1\right) ^{p\beta /2}}du  \notag \\
&=&2\xi ^{p\left( n-\alpha \beta \right) }\left[ \int_{0}^{\infty }\left(
1+u\right) ^{rp+1}\,\frac{u^{np-1}}{\left( u^{2\alpha }+1\right) ^{p\beta /2}%
}du-\int_{0}^{\infty }\,\frac{u^{np-1}}{\left( u^{2\alpha }+1\right)
^{p\beta /2}}du\right] .
\end{eqnarray}%
Thus by (17) and (29) we obtain 
\begin{equation}
\mathcal{J}=2\xi ^{p\left( n-\alpha \beta \right) }\tau .
\end{equation}%
We notice that%
\begin{eqnarray*}
0 &<&\tau <\int_{0}^{\infty }\frac{\left( 1+u\right) ^{rp+1}u^{np-1}}{\left(
u^{2\alpha }+1\right) ^{p\beta /2}}du \\
&<&\int_{0}^{\infty }\frac{\left( 1+u\right) ^{rp+1}\left( 1+u\right) ^{np-1}%
}{\left( u^{2\alpha }+1\right) ^{p\beta /2}}du \\
&=&\int_{0}^{\infty }\frac{\left( 1+u\right) ^{p\left( n+r\right) }}{\left(
u^{2\alpha }+1\right) ^{p\beta /2}}du=:I_{1}.
\end{eqnarray*}%
Also call%
\begin{equation*}
K:=\int_{0}^{1}\frac{\left( 1+u\right) ^{p\left( n+r\right) }}{\left(
u^{2\alpha }+1\right) ^{p\beta /2}}du<\infty .
\end{equation*}%
Then we can write%
\begin{equation*}
I_{1}=K+\int_{1}^{\infty }\frac{\left( 1+u\right) ^{p\left( n+r\right) }}{%
\left( u^{2\alpha }+1\right) ^{p\beta /2}}du<K+2^{p\left( n+r\right)
}\int_{1}^{\infty }\frac{u^{p\left( n+r\right) }}{\left( u^{2\alpha
}+1\right) ^{p\beta /2}}du=K+2^{p\left( n+r\right) }I_{2},
\end{equation*}%
where $I_{2}:=\int_{1}^{\infty }\frac{u^{p\left( n+r\right) }}{\left(
u^{2\alpha }+1\right) ^{p\beta /2}}du.$

Since $\frac{1}{1+u^{2\alpha }}<\frac{1}{u^{2\alpha }},$ we have $\frac{1}{%
\left( 1+u^{2\alpha }\right) ^{p\beta /2}}<\frac{1}{u^{p\alpha \beta }},$
for $u\in \lbrack 1,\infty ).$

So we get 
\begin{eqnarray*}
I_{2} &<&\int_{1}^{\infty }u^{p\left( n+r-\alpha \beta \right)
}du=\lim_{\varepsilon \rightarrow \infty }\int_{1}^{\varepsilon }u^{p\left(
n+r-\alpha \beta \right) }du \\
&=&\lim_{\varepsilon \rightarrow \infty }\left( \frac{\varepsilon ^{p\left(
n+r-\alpha \beta \right) +1}-1}{p\left( n+r-\alpha \beta \right) +1}\right) =%
\frac{-1}{p\left( n+r-\alpha \beta \right) +1},
\end{eqnarray*}%
which is a positive number since $\beta >\frac{1}{\alpha }\left( \frac{1}{p}%
+n+r\right) .$

Consequently $I_{2}$ is finite, so is $I_{1},$ proving $\tau <\infty .$

Using (27) and (30) we get 
\begin{eqnarray}
(\ast \ast ) &=&\left( \frac{\xi c_{1}}{rp+1}\right) \left( \omega
_{r}(f^{(n)},\xi )_{p}\right) ^{p}2\xi ^{p\left( n-\alpha \beta \right) }\tau
\notag \\
&=&\frac{2\alpha \left[ \Gamma \left( \beta \right) \right] ^{p}\Gamma
\left( \frac{q\beta }{2}-\frac{1}{2\alpha }\right) ^{\frac{p}{q}}\tau }{%
\left( rp+1\right) \Gamma \left( \frac{1}{2\alpha }\right) \Gamma \left(
\beta -\frac{1}{2\alpha }\right) ^{p}\Gamma \left( \frac{q\beta }{2}\right)
^{\frac{p}{q}}((n-1)!)^{p}(q(n-1)+1)^{p/q}}\xi ^{pn}\left( \omega
_{r}(f^{(n)},\xi )_{p}\right) ^{p}.
\end{eqnarray}%
I.e. we have established that 
\begin{equation}
I\leq \frac{2\alpha \left[ \Gamma \left( \beta \right) \right] ^{p}\Gamma
\left( \frac{q\beta }{2}-\frac{1}{2\alpha }\right) ^{\frac{p}{q}}\tau }{%
\left( rp+1\right) \Gamma \left( \frac{1}{2\alpha }\right) \Gamma \left(
\beta -\frac{1}{2\alpha }\right) ^{p}\Gamma \left( \frac{q\beta }{2}\right)
^{\frac{p}{q}}((n-1)!)^{p}(q(n-1)+1)^{p/q}}\xi ^{pn}\left( \omega
_{r}(f^{(n)},\xi )_{p}\right) ^{p}\,.
\end{equation}%
That is finishing the proof of the theorem. \hfill $\blacksquare $

\medskip The counterpart of Theorem 1 follows, case of $p=1.$

\medskip \noindent \textbf{Theorem 2.} \textit{Let }$f\in C^{n}(\mathbb{R})$
and $f^{(n)}\in L_{1}(\mathbb{R}),n\in \mathbb{N},$ $\alpha \in \mathbb{N},$ 
$\beta >\frac{n+r+1}{2\alpha }$.\textit{\ Then} 
\begin{eqnarray}
\Vert \Delta (x)\Vert _{1} &\leq &\frac{1}{\left( r+1\right) (n-1)!\Gamma
\left( \frac{1}{2\alpha }\right) \Gamma \left( \beta -\frac{1}{2\alpha }%
\right) } \\
&&\cdot \left[ \sum_{k=1}^{r+1}{\binom{{r+1}}{k}}\Gamma \left( \frac{n+k}{%
2\alpha }\right) \Gamma \left( \beta -\frac{n+k}{2\alpha }\right) \right]
\omega _{r}(f^{(n)},\xi )_{1}\xi ^{n}.  \notag
\end{eqnarray}%
\textit{Hence as }$\xi \rightarrow 0$\textit{\ we obtain} $\Vert \Delta
(x)\Vert _{1}\rightarrow 0$.

If additionally $f^{(2m)}\in L_{1}(\mathbb{R}),m=1,2,\ldots ,\left\lfloor 
\frac{n}{2}\right\rfloor $ then $\left\Vert M_{r,\xi }(f)-f\right\Vert
_{1}\rightarrow 0,$ as $\xi \rightarrow 0.$

\smallskip \noindent \textbf{Proof.} It follows 
\begin{eqnarray}
|\Delta (x)| &=&W\left\vert \int_{-\infty }^{\infty }\mathcal{R}_{n}(0,t,x)%
\frac{1}{\left( t^{2\alpha }+\xi ^{2\alpha }\right) ^{\beta }}dt\right\vert 
\notag \\
&\leq &W\int_{-\infty }^{\infty }\left\vert \mathcal{R}_{n}(0,t,x)\right%
\vert \frac{1}{\left( t^{2\alpha }+\xi ^{2\alpha }\right) ^{\beta }}dt 
\notag \\
&\leq &W\int_{-\infty }^{\infty }\left( \int_{0}^{|t|}\frac{(|t|-w)^{n-1}}{%
(n-1)!}|\tau (sign(t)\cdot w,x)|dw\right) \frac{1}{\left( t^{2\alpha }+\xi
^{2\alpha }\right) ^{\beta }}dt.
\end{eqnarray}%
Thus 
\begin{gather}
\Vert \Delta (x)\Vert _{1}=\int_{-\infty }^{\infty }|\Delta (x)|dx\leq
W\cdot \text{ \ \ \ \ \ \ \ \ \ \ \ \ \ \ \ \ \ \ \ \ \ \ \ \ \ \ \ \ \ \ \
\ \ \ \ \ \ \ \ \ \ \ \ \ \ \ \ \ \ \ \ \ \ \ \ \ \ \ \ \ \ \ \ \ \ \ \ } \\
\int_{-\infty }^{\infty }\left( \int_{-\infty }^{\infty }\left(
\int_{0}^{|t|}\frac{(|t|-w)^{n-1}}{(n-1)!}|\tau (sign(t)\cdot w,x)|dw\right) 
\frac{1}{\left( t^{2\alpha }+\xi ^{2\alpha }\right) ^{\beta }}dt\right) dx 
\notag \\
=:(\ast )\text{ \ \ \ \ \ \ \ \ \ \ \ \ \ \ \ \ \ \ \ \ \ \ \ \ \ \ \ \ \ \
\ \ \ \ \ \ \ \ \ \ \ \ \ \ \ \ \ \ \ \ \ \ \ \ \ \ \ \ \ \ \ \ \ \ \ \ \ \
\ \ \ \ \ \ \ \ \ \ \ \ \ \ \ \ \ \ \ \ \ \ \ \ \ \ \ \ \ \ \ \ \ }  \notag
\end{gather}%
But we see that 
\begin{equation}
\int_{0}^{|t|}\frac{(|t|-w)^{n-1}}{(n-1)!}|\tau (sign(t)\cdot w,x)|dw\leq 
\frac{|t|^{n-1}}{(n-1)!}\int_{0}^{|t|}|\tau (sign(t)\cdot w,x)|dw.
\end{equation}%
Therefore it holds 
\begin{eqnarray}
(\ast ) &\leq &W\int_{-\infty }^{\infty }\left( \int_{-\infty }^{\infty
}\left( \frac{|t|^{n-1}}{(n-1)!}\int_{0}^{|t|}|\tau (sign(t)\cdot
w,x)|dw\right) \frac{1}{\left( t^{2\alpha }+\xi ^{2\alpha }\right) ^{\beta }}%
dt\right) dx  \notag \\
&=&\frac{W}{(n-1)!}\left( \int_{-\infty }^{\infty }\left(
\int_{0}^{|t|}\left( \int_{-\infty }^{\infty }|\tau (sign(t)\cdot
w,x)|dx\right) dw\right) \frac{|t|^{n-1}}{\left( t^{2\alpha }+\xi ^{2\alpha
}\right) ^{\beta }}dt\right)  \notag \\
&\leq &\frac{W}{(n-1)!}\left( \int_{-\infty }^{\infty }\left(
\int_{0}^{|t|}\omega _{r}(f^{(n)},w)_{1}dw\right) \frac{|t|^{n-1}}{\left(
t^{2\alpha }+\xi ^{2\alpha }\right) ^{\beta }}dt\right) .
\end{eqnarray}%
I.e. we get 
\begin{equation}
\Vert \Delta (x)\Vert _{1}\leq \frac{W}{(n-1)!}\left( \int_{-\infty
}^{\infty }\left( \int_{0}^{|t|}\omega _{r}(f^{(n)},w)_{1}dw\right) \frac{%
|t|^{n-1}}{\left( t^{2\alpha }+\xi ^{2\alpha }\right) ^{\beta }}dt\right) .
\end{equation}

Consequently we have 
\begin{eqnarray}
\Vert \Delta (x)\Vert _{1} &\leq &\frac{W\omega _{r}(f^{(n)},\xi )_{1}}{%
(n-1)!}\left( \int_{-\infty }^{\infty }\left( \int_{0}^{|t|}\left( 1+\frac{w%
}{\xi }\right) ^{r}dw\right) \frac{|t|^{n-1}}{\left( t^{2\alpha }+\xi
^{2\alpha }\right) ^{\beta }}dt\right)  \notag \\
&=&\frac{2\xi W\omega _{r}(f^{(n)},\xi )_{1}}{\left( r+1\right) (n-1)!}%
\left( \int_{0}^{\infty }\left( \left( 1+\frac{t}{\xi }\right)
^{r+1}-1\right) \frac{t^{n-1}}{\left( t^{2\alpha }+\xi ^{2\alpha }\right)
^{\beta }}dt\right)  \notag \\
&=&\frac{2\Gamma \left( \beta \right) \alpha \xi ^{2\alpha \beta }\omega
_{r}(f^{(n)},\xi )_{1}}{\left( r+1\right) (n-1)!\Gamma \left( \frac{1}{%
2\alpha }\right) \Gamma \left( \beta -\frac{1}{2\alpha }\right) }\left(
\int_{0}^{\infty }\left( \left( 1+\frac{t}{\xi }\right) ^{r+1}-1\right) 
\frac{t^{n-1}}{\left( t^{2\alpha }+\xi ^{2\alpha }\right) ^{\beta }}%
dt\right) .
\end{eqnarray}%
We have gotten so far 
\begin{equation}
\Vert \Delta (x)\Vert _{1}\leq \frac{2\Gamma \left( \beta \right) \alpha \xi
^{2\alpha \beta }\omega _{r}(f^{(n)},\xi )_{1}\cdot \lambda }{\left(
r+1\right) (n-1)!\Gamma \left( \frac{1}{2\alpha }\right) \Gamma \left( \beta
-\frac{1}{2\alpha }\right) }\,,
\end{equation}%
where 
\begin{equation}
\lambda :=\int_{0}^{\infty }\left( \left( 1+\frac{t}{\xi }\right)
^{r+1}-1\right) \frac{t^{n-1}}{\left( t^{2\alpha }+\xi ^{2\alpha }\right)
^{\beta }}dt.
\end{equation}

One easily finds that 
\begin{eqnarray}
\lambda &=&\int_{0}^{\infty }\left( \sum_{k=1}^{r+1}{\binom{{r+1}}{k}}\left( 
\frac{t}{\xi }\right) ^{k}\right) \frac{t^{n-1}}{\left( t^{2\alpha }+\xi
^{2\alpha }\right) ^{\beta }}dt  \notag \\
&=&\xi ^{n-2\alpha \beta }\sum_{k=1}^{r+1}{\binom{{r+1}}{k}}\int_{0}^{\infty
}\frac{T^{n+k-1}}{\left( T^{2\alpha }+1\right) ^{\beta }}dT  \notag \\
&=&\xi ^{n-2\alpha \beta }\sum_{k=1}^{r+1}{\binom{{r+1}}{k}}K_{n+k}.
\end{eqnarray}%
Where 
\begin{equation}
K_{n+k}:=\int_{0}^{\infty }\frac{T^{n+k-1}}{\left( T^{2\alpha }+1\right)
^{\beta }}dT=\frac{\Gamma \left( \frac{n+k}{2\alpha }\right) \Gamma \left(
\beta -\frac{n+k}{2\alpha }\right) }{\Gamma \left( \beta \right) 2\alpha }.
\end{equation}%
\begin{equation*}
\Vert \Delta (x)\Vert _{1}\leq \frac{1}{\left( r+1\right) (n-1)!\Gamma
\left( \frac{1}{2\alpha }\right) \Gamma \left( \beta -\frac{1}{2\alpha }%
\right) }\left[ \sum_{k=1}^{r+1}{\binom{{r+1}}{k}}\Gamma \left( \frac{n+k}{%
2\alpha }\right) \Gamma \left( \beta -\frac{n+k}{2\alpha }\right) \right]
\omega _{r}(f^{(n)},\xi )_{1}\xi ^{n}.
\end{equation*}

We have proved (33).\hfill $\blacksquare $

\medskip The case $n=0$ is met next.

\medskip \noindent \textbf{Proposition 1.} \textit{Let }$p,q>1$\textit{\
such that }$\frac{1}{p}+\frac{1}{q}=1,$ $\alpha \in \mathbb{N},$ $\beta >%
\frac{1}{\alpha }\left( r+\frac{1}{p}\right) $\textit{\ and the rest as
above. Then} 
\begin{equation}
\Vert M_{r,\xi }(f)-f\Vert _{p}\leq \frac{\left( 2\alpha \right) ^{\frac{1}{p%
}}\left[ \Gamma \left( \beta \right) \right] \Gamma \left( \frac{q\beta }{2}-%
\frac{1}{2\alpha }\right) ^{\frac{1}{q}}\theta ^{\frac{1}{p}}}{\Gamma \left( 
\frac{1}{2\alpha }\right) ^{\frac{1}{p}}\Gamma \left( \beta -\frac{1}{%
2\alpha }\right) \Gamma \left( \frac{q\beta }{2}\right) ^{\frac{1}{q}}}%
\omega _{r}(f,\xi )_{p},
\end{equation}%
\textit{where} 
\begin{equation}
0<\theta :=\int_{0}^{\infty }\left( 1+t\right) ^{rp}\frac{1}{\left(
t^{2\alpha }+1\right) ^{p\beta /2}}dt<\infty .
\end{equation}%
Hence as $\xi \rightarrow 0$ we obtain $M_{r,\xi }$ $\rightarrow $ unit
operator $I$ in the $L_{p}$ norm, $p>1$.

\smallskip \noindent \textbf{Proof.} By (3) we notice that,%
\begin{align}
M_{r,\xi }(f;x)-f(x)& =W\left( \sum_{j=0}^{r}\alpha _{j}\int_{-\infty
}^{\infty }(f(x+jt)-f(x))\frac{1}{\left( t^{2\alpha }+\xi ^{2\alpha }\right)
^{\beta }}dt\right)  \notag \\
& =W\left( \int_{-\infty }^{\infty }\left( \sum_{j=0}^{r}\alpha _{j}\left(
f(x+jt)-f(x)\right) \right) \frac{1}{\left( t^{2\alpha }+\xi ^{2\alpha
}\right) ^{\beta }}dt\right)  \notag \\
& =W\left( \int_{-\infty }^{\infty }\left( \sum_{j=1}^{r}\alpha
_{j}f(x+jt)-\sum_{j=1}^{r}\alpha _{j}f(x)\right) \frac{1}{\left( t^{2\alpha
}+\xi ^{2\alpha }\right) ^{\beta }}dt\right)  \notag
\end{align}%
\begin{align}
& =W\left( \int_{-\infty }^{\infty }\left( \sum_{j=1}^{r}\left( (-1)^{r-j}{%
\binom{r}{j}}j^{-n}\right) f(x+jt)-\sum_{j=1}^{r}\left( (-1)^{r-j}{\binom{r}{%
j}}j^{-n}\right) f(x)\right) \frac{1}{\left( t^{2\alpha }+\xi ^{2\alpha
}\right) ^{\beta }}dt\right)  \notag \\
& =W\left( \int_{-\infty }^{\infty }\left( \sum_{j=1}^{r}\left( (-1)^{r-j}{%
\binom{r}{j}}\right) f(x+jt)-\sum_{j=1}^{r}\left( (-1)^{r-j}{\binom{r}{j}}%
\right) f(x)\right) \frac{1}{\left( t^{2\alpha }+\xi ^{2\alpha }\right)
^{\beta }}dt\right)  \notag \\
& \overset{(9)}{=}W\left( \int_{-\infty }^{\infty }\left(
\sum_{j=1}^{r}\left( (-1)^{r-j}{\binom{r}{j}}\right) f(x+jt)+\left(
(-1)^{r-0}{\binom{r}{0}}\right) f(x+0t)\right) \frac{1}{\left( t^{2\alpha
}+\xi ^{2\alpha }\right) ^{\beta }}dt\right)  \notag \\
& =W\left( \int_{-\infty }^{\infty }\left( \sum_{j=0}^{r}(-1)^{r-j}{\binom{r%
}{j}}f(x+jt)\right) \frac{1}{\left( t^{2\alpha }+\xi ^{2\alpha }\right)
^{\beta }}dt\right)  \notag \\
& \overset{(6)}{=}W\left( \int_{-\infty }^{\infty }\left( \Delta
_{t}^{r}f\right) (x)\frac{1}{\left( t^{2\alpha }+\xi ^{2\alpha }\right)
^{\beta }}dt\right) .
\end{align}%
And then 
\begin{equation}
\left\vert M_{r,\xi }(f;x)-f(x)\right\vert \leq W\left( \int_{-\infty
}^{\infty }\left\vert \Delta _{t}^{r}f(x)\right\vert \frac{1}{\left(
t^{2\alpha }+\xi ^{2\alpha }\right) ^{\beta }}dt\right) .
\end{equation}%
We next estimate%
\begin{eqnarray*}
&&\int_{-\infty }^{\infty }|M_{r,\xi }(f;x)-f(x)|^{p}dx\leq \int_{-\infty
}^{\infty }\left( W\right) ^{p}\left( \int_{-\infty }^{\infty }\left\vert
\Delta _{t}^{r}f(x)\right\vert \frac{1}{\left( t^{2\alpha }+\xi ^{2\alpha
}\right) ^{\beta }}dt\right) ^{p}dx \\
&=&\left( W\right) ^{p}\int_{-\infty }^{\infty }\left( \int_{-\infty
}^{\infty }\left( \left\vert \Delta _{t}^{r}f(x)\right\vert \frac{1}{\left(
t^{2\alpha }+\xi ^{2\alpha }\right) ^{\beta /2}}\right) \left( \frac{1}{%
\left( t^{2\alpha }+\xi ^{2\alpha }\right) ^{\beta /2}}\right) dt\right)
^{p}dx \\
&\leq &\left( W\right) ^{p}\int_{-\infty }^{\infty }\left( \left(
\int_{-\infty }^{\infty }\left( \left\vert \Delta _{t}^{r}f(x)\right\vert 
\frac{1}{\left( t^{2\alpha }+\xi ^{2\alpha }\right) ^{\beta /2}}\right)
^{p}dt\right) ^{\frac{1}{p}}\left( \int_{-\infty }^{\infty }\left( \frac{1}{%
\left( t^{2\alpha }+\xi ^{2\alpha }\right) ^{\beta /2}}\right) ^{q}dt\right)
^{\frac{1}{q}}\right) ^{p}dx
\end{eqnarray*}%
\begin{eqnarray*}
&=&\left( W\right) ^{p}\int_{-\infty }^{\infty }\left( \int_{-\infty
}^{\infty }\left\vert \Delta _{t}^{r}f(x)\right\vert ^{p}\frac{1}{\left(
t^{2\alpha }+\xi ^{2\alpha }\right) ^{p\beta /2}}dt\right) \left(
\int_{-\infty }^{\infty }\frac{1}{\left( t^{2\alpha }+\xi ^{2\alpha }\right)
^{q\beta /2}}dt\right) ^{\frac{p}{q}}dx \\
&=&\left( W\right) ^{p}\int_{-\infty }^{\infty }\left( \int_{-\infty
}^{\infty }\left\vert \Delta _{t}^{r}f(x)\right\vert ^{p}\frac{1}{\left(
t^{2\alpha }+\xi ^{2\alpha }\right) ^{p\beta /2}}dt\right) \left( \frac{%
\Gamma \left( \frac{1}{2\alpha }\right) \Gamma \left( \frac{q\beta }{2}-%
\frac{1}{2\alpha }\right) }{\Gamma \left( \frac{q\beta }{2}\right) \alpha
\xi ^{q\alpha \beta -1}}\right) ^{\frac{p}{q}}dx \\
&=&\left( W\right) ^{p}\left( \frac{\Gamma \left( \frac{1}{2\alpha }\right)
\Gamma \left( \frac{q\beta }{2}-\frac{1}{2\alpha }\right) }{\Gamma \left( 
\frac{q\beta }{2}\right) \alpha \xi ^{q\alpha \beta -1}}\right) ^{\frac{p}{q}%
}\int_{-\infty }^{\infty }\left( \int_{-\infty }^{\infty }\left\vert \Delta
_{t}^{r}f(x)\right\vert ^{p}\frac{1}{\left( t^{2\alpha }+\xi ^{2\alpha
}\right) ^{p\beta /2}}dt\right) dx
\end{eqnarray*}%
\begin{eqnarray}
&=&\frac{\alpha \xi ^{\alpha \beta p-1}\left[ \Gamma \left( \beta \right) %
\right] ^{p}\Gamma \left( \frac{q\beta }{2}-\frac{1}{2\alpha }\right) ^{%
\frac{p}{q}}}{\Gamma \left( \frac{1}{2\alpha }\right) \Gamma \left( \beta -%
\frac{1}{2\alpha }\right) ^{p}\Gamma \left( \frac{q\beta }{2}\right) ^{\frac{%
p}{q}}}\int_{-\infty }^{\infty }\left( \int_{-\infty }^{\infty }\left\vert
\Delta _{t}^{r}f(x)\right\vert ^{p}\frac{1}{\left( t^{2\alpha }+\xi
^{2\alpha }\right) ^{p\beta /2}}dx\right) dt  \notag \\
&=&\frac{\alpha \xi ^{\alpha \beta p-1}\left[ \Gamma \left( \beta \right) %
\right] ^{p}\Gamma \left( \frac{q\beta }{2}-\frac{1}{2\alpha }\right) ^{%
\frac{p}{q}}}{\Gamma \left( \frac{1}{2\alpha }\right) \Gamma \left( \beta -%
\frac{1}{2\alpha }\right) ^{p}\Gamma \left( \frac{q\beta }{2}\right) ^{\frac{%
p}{q}}}\int_{-\infty }^{\infty }\left( \int_{-\infty }^{\infty }\left\vert
\Delta _{t}^{r}f(x)\right\vert ^{p}dx\right) \frac{1}{\left( t^{2\alpha
}+\xi ^{2\alpha }\right) ^{p\beta /2}}dt  \notag \\
&\leq &\frac{\alpha \xi ^{\alpha \beta p-1}\left[ \Gamma \left( \beta
\right) \right] ^{p}\Gamma \left( \frac{q\beta }{2}-\frac{1}{2\alpha }%
\right) ^{\frac{p}{q}}}{\Gamma \left( \frac{1}{2\alpha }\right) \Gamma
\left( \beta -\frac{1}{2\alpha }\right) ^{p}\Gamma \left( \frac{q\beta }{2}%
\right) ^{\frac{p}{q}}}\int_{-\infty }^{\infty }\omega _{r}(f,\left\vert
t\right\vert )_{p}^{p}\frac{1}{\left( t^{2\alpha }+\xi ^{2\alpha }\right)
^{p\beta /2}}dt  \notag
\end{eqnarray}%
\begin{eqnarray}
&=&\frac{2\alpha \xi ^{\alpha \beta p-1}\left[ \Gamma \left( \beta \right) %
\right] ^{p}\Gamma \left( \frac{q\beta }{2}-\frac{1}{2\alpha }\right) ^{%
\frac{p}{q}}}{\Gamma \left( \frac{1}{2\alpha }\right) \Gamma \left( \beta -%
\frac{1}{2\alpha }\right) ^{p}\Gamma \left( \frac{q\beta }{2}\right) ^{\frac{%
p}{q}}}\int_{0}^{\infty }\omega _{r}(f,t)_{p}^{p}\frac{1}{\left( t^{2\alpha
}+\xi ^{2\alpha }\right) ^{p\beta /2}}dt  \notag \\
&\leq &\frac{2\alpha \xi ^{\alpha \beta p-1}\left[ \Gamma \left( \beta
\right) \right] ^{p}\Gamma \left( \frac{q\beta }{2}-\frac{1}{2\alpha }%
\right) ^{\frac{p}{q}}}{\Gamma \left( \frac{1}{2\alpha }\right) \Gamma
\left( \beta -\frac{1}{2\alpha }\right) ^{p}\Gamma \left( \frac{q\beta }{2}%
\right) ^{\frac{p}{q}}}\omega _{r}(f,\xi )_{p}^{p}\int_{0}^{\infty }\left( 1+%
\frac{t}{\xi }\right) ^{rp}\frac{1}{\left( t^{2\alpha }+\xi ^{2\alpha
}\right) ^{p\beta /2}}dt  \notag \\
&=&\frac{2\alpha \xi ^{\alpha \beta p-1}\left[ \Gamma \left( \beta \right) %
\right] ^{p}\Gamma \left( \frac{q\beta }{2}-\frac{1}{2\alpha }\right) ^{%
\frac{p}{q}}}{\Gamma \left( \frac{1}{2\alpha }\right) \Gamma \left( \beta -%
\frac{1}{2\alpha }\right) ^{p}\Gamma \left( \frac{q\beta }{2}\right) ^{\frac{%
p}{q}}}\omega _{r}(f,\xi )_{p}^{p}\int_{0}^{\infty }\left( 1+T\right) ^{rp}%
\frac{1}{\left( T^{2\alpha }+1\right) ^{p\beta /2}\xi ^{\alpha p\beta }}\xi
dT  \notag \\
&=&\frac{2\alpha \left[ \Gamma \left( \beta \right) \right] ^{p}\Gamma
\left( \frac{q\beta }{2}-\frac{1}{2\alpha }\right) ^{\frac{p}{q}}}{\Gamma
\left( \frac{1}{2\alpha }\right) \Gamma \left( \beta -\frac{1}{2\alpha }%
\right) ^{p}\Gamma \left( \frac{q\beta }{2}\right) ^{\frac{p}{q}}}\omega
_{r}(f,\xi )_{p}^{p}\int_{0}^{\infty }\left( 1+t\right) ^{rp}\frac{1}{\left(
t^{2\alpha }+1\right) ^{p\beta /2}}dt.
\end{eqnarray}%
We have established (44).

We also notice that%
\begin{eqnarray*}
\theta &=&\int_{0}^{\infty }\frac{\left( 1+t\right) ^{rp}}{\left(
1+t^{2\alpha }\right) ^{p\beta /2}}dt=\int_{0}^{1}\frac{\left( 1+t\right)
^{rp}}{\left( 1+t^{2\alpha }\right) ^{p\beta /2}}dt+\int_{1}^{\infty }\frac{%
\left( 1+t\right) ^{rp}}{\left( 1+t^{2\alpha }\right) ^{p\beta /2}}dt \\
&<&\int_{0}^{1}\frac{\left( 1+t\right) ^{rp}}{\left( 1+t^{2\alpha }\right)
^{p\beta /2}}dt+2^{rp}\int_{1}^{\infty }\frac{t^{rp}}{\left( 1+t^{2\alpha
}\right) ^{p\beta /2}}dt \\
&<&\int_{0}^{1}\frac{\left( 1+t\right) ^{rp}}{\left( 1+t^{2\alpha }\right)
^{p\beta /2}}dt+2^{rp}\int_{1}^{\infty }t^{p\left( r-\alpha \beta \right) }dt
\\
&=&\int_{0}^{1}\frac{\left( 1+t\right) ^{rp}}{\left( 1+t^{2\alpha }\right)
^{p\beta /2}}dt-\frac{2^{rp}}{p\left( r-\alpha \beta \right) +1},
\end{eqnarray*}%
the last, since $\beta >\frac{1}{\alpha }\left( r+\frac{1}{p}\right) ,$ is a
finite positive constant. Thus $0<\theta <$ $\infty .$\hfill $\blacksquare $

\medskip We also give

\medskip \noindent \textbf{Proposition 2.} \textit{Assume }$\beta >\frac{r+1%
}{2\alpha }.$ It holds 
\begin{equation}
\Vert M_{r,\xi }f-f\Vert _{1}\leq \frac{2\alpha \Gamma \left( \beta \right) 
}{\Gamma \left( \frac{1}{2\alpha }\right) \Gamma \left( \beta -\frac{1}{%
2\alpha }\right) }\left( \int_{0}^{\infty }\left( 1+t\right) ^{r}\frac{1}{%
\left( t^{2\alpha }+1\right) ^{\beta }}dt\right) \omega _{r}(f,\xi )_{1}.
\end{equation}%
\textit{Hence as $\xi \rightarrow 0$ we get $M_{r,\xi }\rightarrow I$ in the 
$L_{1}$ norm.}

\smallskip \noindent \textbf{Proof.} By (47) we have again 
\begin{equation*}
\left\vert M_{r,\xi }(f;x)-f(x)\right\vert \leq W\left( \int_{-\infty
}^{\infty }\left\vert \Delta _{t}^{r}f(x)\right\vert \frac{1}{\left(
t^{2\alpha }+\xi ^{2\alpha }\right) ^{\beta }}dt\right) .
\end{equation*}%
Next we estimate 
\begin{eqnarray}
\int_{-\infty }^{\infty }\left\vert M_{r,\xi }(f;x)-f(x)\right\vert dx &\leq
&\;W\int_{-\infty }^{\infty }\left( \int_{-\infty }^{\infty }\left\vert
\Delta _{t}^{r}f(x)\right\vert \frac{1}{\left( t^{2\alpha }+\xi ^{2\alpha
}\right) ^{\beta }}dt\right) dx  \notag \\
&=&W\int_{-\infty }^{\infty }\left( \int_{-\infty }^{\infty }\left\vert
\Delta _{t}^{r}f(x)\right\vert dx\right) \frac{1}{\left( t^{2\alpha }+\xi
^{2\alpha }\right) ^{\beta }}dt  \notag \\
&\leq &W\int_{-\infty }^{\infty }\omega _{r}(f,|t|)_{1}\frac{1}{\left(
t^{2\alpha }+\xi ^{2\alpha }\right) ^{\beta }}dt  \notag \\
&\leq &\;W2\omega _{r}(f,\xi )_{1}\int_{0}^{\infty }\left( 1+\frac{t}{\xi }%
\right) ^{r}\frac{1}{\left( t^{2\alpha }+\xi ^{2\alpha }\right) ^{\beta }}dt
\notag \\
&=&\;\;\frac{\Gamma \left( \beta \right) \xi ^{2\alpha \beta -1}2\alpha }{%
\Gamma \left( \frac{1}{2\alpha }\right) \Gamma \left( \beta -\frac{1}{%
2\alpha }\right) }\omega _{r}(f,\xi )_{1}\int_{0}^{\infty }\xi \left(
1+t\right) ^{r}\frac{1}{\left( t^{2\alpha }+1\right) ^{\beta }\xi ^{2\alpha
\beta }}dt  \notag \\
&=&\;\frac{\Gamma \left( \beta \right) 2\alpha }{\Gamma \left( \frac{1}{%
2\alpha }\right) \Gamma \left( \beta -\frac{1}{2\alpha }\right) }\omega
_{r}(f,\xi )_{1}\int_{0}^{\infty }\left( 1+t\right) ^{r}\frac{1}{\left(
t^{2\alpha }+1\right) ^{\beta }}dt.
\end{eqnarray}%
We have proved (49).

We also notice that 
\begin{eqnarray*}
0 &<&\int_{0}^{\infty }\left( 1+t\right) ^{r}\frac{1}{\left( t^{2\alpha
}+1\right) ^{\beta }}dt \\
&=&\int_{0}^{1}\frac{\left( 1+t\right) ^{r}}{\left( t^{2\alpha }+1\right)
^{\beta }}dt+\int_{1}^{\infty }\frac{\left( 1+t\right) ^{r}}{\left(
t^{2\alpha }+1\right) ^{\beta }}dt \\
&<&\int_{0}^{1}\frac{\left( 1+t\right) ^{r}}{\left( t^{2\alpha }+1\right)
^{\beta }}dt+2^{r}\int_{1}^{\infty }t^{r-2\alpha \beta }dt \\
&=&\int_{0}^{1}\frac{\left( 1+t\right) ^{r}}{\left( t^{2\alpha }+1\right)
^{\beta }}dt-\frac{2^{r}}{\left( r-2\alpha \beta +1\right) },
\end{eqnarray*}%
which is a positive finite constant. \hfill $\blacksquare $

\medskip In the next we consider $f\in C^{n}(\mathbb{R})$ and $f^{(n)}\in
L_{p}(\mathbb{R}),n=0$ or $n\geq 2$ even, $1\leq p<\infty $ and the similar 
\textit{smooth singular operator of symmetric convolution type} 
\begin{equation}
M_{\xi }(f;x)=W\int_{-\infty }^{\infty }f(x+y)\frac{1}{\left( y^{2\alpha
}+\xi ^{2\alpha }\right) ^{\beta }}dy,\ \ \hbox{for all }x\in \mathbb{R},\
\xi >0.
\end{equation}%
That is 
\begin{equation*}
M_{\xi }(f;x)=W\int_{0}^{\infty }\left( f(x+y)+f(x-y)\right) \frac{1}{\left(
y^{2\alpha }+\xi ^{2\alpha }\right) ^{\beta }}dy,
\end{equation*}%
for all $x\in \mathbb{R}$, $\xi >0$. Notice that $M_{1,\xi }=M_{\xi }$. Let
the central second order difference 
\begin{equation}
(\tilde{\Delta}_{y}^{2}f)(x):=f(x+y)+f(x-y)-2f(x).
\end{equation}%
Notice that 
\begin{equation*}
(\tilde{\Delta}_{-y}^{2}f)(x)=(\tilde{\Delta}_{y}^{2}f)(x).
\end{equation*}%
When $n\geq 2$ even using Taylor's formula with Cauchy remainder we
eventually find 
\begin{equation}
(\tilde{\Delta}_{y}^{2}f)(x)=2\sum_{\rho =1}^{n/2}\frac{f^{(2\rho )}(x)}{%
(2\rho )!}y^{2\rho }+\mathcal{R}_{1}(x),
\end{equation}%
where 
\begin{equation}
\mathcal{R}_{1}(x):=\int_{0}^{y}(\tilde{\Delta}_{t}^{2}f^{(n)})(x)\frac{%
(y-t)^{n-1}}{(n-1)!}dt.
\end{equation}%
Notice that 
\begin{equation}
M_{\xi }(f;x)-f(x)=W\int_{0}^{\infty }(\tilde{\Delta}_{y}^{2}f(x))\frac{1}{%
\left( y^{2\alpha }+\xi ^{2\alpha }\right) ^{\beta }}dy.
\end{equation}

Furthermore by (4), (53) and (55) we easily see that 
\begin{eqnarray}
K(x) &:&=M_{\xi }(f;x)-f(x)-\sum_{\rho =1}^{n/2}\frac{f^{(2\rho )}(x)}{%
(2\rho )!}\frac{\Gamma \left( \frac{2\rho +1}{2\alpha }\right) \Gamma \left(
\beta -\frac{2\rho +1}{2\alpha }\right) }{\Gamma \left( \frac{1}{2\alpha }%
\right) \Gamma \left( \beta -\frac{1}{2\alpha }\right) }\xi ^{2\rho } \\
&=&W\int_{0}^{\infty }\left[ \int_{0}^{y}(\tilde{\Delta}_{t}^{2}f^{(n)})(x)%
\frac{(y-t)^{n-1}}{(n-1)!}dt\right] \frac{1}{\left( y^{2\alpha }+\xi
^{2\alpha }\right) ^{\beta }}dy,  \notag
\end{eqnarray}%
where $\beta >\frac{(n+1)}{2\alpha }.$

Therefore we have 
\begin{equation}
|K(x)|\leq W\int_{0}^{\infty }\left( \int_{0}^{y}\left\vert \tilde{\Delta}%
_{t}^{2}f^{(n)}\right\vert (x)\frac{(y-t)^{n-1}}{(n-1)!}dt\right) \frac{1}{%
\left( y^{2\alpha }+\xi ^{2\alpha }\right) ^{\beta }}dy.
\end{equation}%
Here we estimate in $L_{p}$ norm, $p\geq 1$, the error function $K(x)$.
Notice that we have $\omega _{2}(f^{(n)},h)_{p}<\infty $, $h>0$, $n=0$ or $%
n\geq 2$ even. Operators $M_{\xi }$ are positive operators.

The related main $L_p$ result here comes next.

\medskip \noindent \textbf{Theorem 3.} \textit{Let }$p,q>1$\textit{\ such
that }$\frac{1}{p}+\frac{1}{q}=1$, $n\geq 2$\textit{\ even,} $\alpha \in 
\mathbb{N},$ $\beta >\frac{1}{\alpha }\left( \frac{1}{p}+n+2\right) $\textit{%
\ and the rest as above. Then} 
\begin{equation}
\Vert K(x)\Vert _{p}\leq \frac{\tilde{\tau}^{1/p}\alpha ^{1/p}\Gamma \left( 
\frac{q\beta }{2}-\frac{1}{2\alpha }\right) ^{1/q}}{2^{\frac{1}{q}}\Gamma
\left( \frac{1}{2\alpha }\right) ^{1/p}\Gamma \left( \beta -\frac{1}{2\alpha 
}\right) \Gamma \left( \frac{q\beta }{2}\right) ^{1/q}(q(n-1)+1)^{1/q}\left(
2p+1\right) ^{1/p}}\frac{\Gamma \left( \beta \right) }{(n-1)!}\xi ^{n}\omega
_{2}(f^{(n)},\xi )_{p},
\end{equation}%
\textit{where} 
\begin{equation}
0<\tilde{\tau}=\int_{0}^{\infty }\left( \left( 1+u\right) ^{2p+1}-1\right)
u^{pn-1}\frac{1}{\left( 1+u^{2\alpha }\right) ^{p\beta /2}}du<\infty .
\end{equation}%
\textit{Hence as }$\xi \rightarrow 0$\textit{\ we get} $\Vert K(x)\Vert
_{p}\rightarrow 0$.

If additionally $f^{(2m)}\in L_{p}(\mathbb{R}),m=1,2,\ldots ,\frac{n}{2}$
then $\left\Vert M_{\xi }(f)-f\right\Vert _{p}\rightarrow 0,$ as $\xi
\rightarrow 0.$

\medskip \noindent \textbf{Proof.} We observe that 
\begin{equation}
|K(x)|^{p}\leq W^{p}\left( \int_{0}^{\infty }\left( \int_{0}^{y}\left\vert 
\tilde{\Delta}_{t}^{2}f^{(n)}\right\vert (x)\frac{(y-t)^{n-1}}{(n-1)!}%
dt\right) \frac{1}{\left( y^{2\alpha }+\xi ^{2\alpha }\right) ^{\beta }}%
dy\right) ^{p}.
\end{equation}%
Call 
\begin{equation}
\tilde{\gamma}(y,x):=\int_{0}^{y}|\tilde{\Delta}_{t}^{2}f^{(n)}(x)|\frac{%
(y-t)^{n-1}}{(n-1)!}dt\geq 0,\text{ \ \ }y\geq 0,
\end{equation}%
then we have 
\begin{equation}
|K(x)|^{p}\leq W^{p}\left( \int_{0}^{\infty }\tilde{\gamma}(y,x)\frac{1}{%
\left( y^{2\alpha }+\xi ^{2\alpha }\right) ^{\beta }}dy\right) ^{p}.
\end{equation}%
Consequently 
\begin{eqnarray}
\Lambda  &:&=\int_{-\infty }^{\infty }|K(x)|^{p}dx\leq W^{p}\int_{-\infty
}^{\infty }\left( \int_{0}^{\infty }\tilde{\gamma}(y,x)\frac{1}{\left(
y^{2\alpha }+\xi ^{2\alpha }\right) ^{\beta }}dy\right) ^{p}dx  \notag \\
&=&W^{p}\int_{-\infty }^{\infty }\left( \int_{0}^{\infty }\tilde{\gamma}(y,x)%
\frac{1}{\left( y^{2\alpha }+\xi ^{2\alpha }\right) ^{\beta /2}}\frac{1}{%
\left( y^{2\alpha }+\xi ^{2\alpha }\right) ^{\beta /2}}dy\right) ^{p}dx 
\notag \\
&\ &\qquad \qquad \hbox{(by H\"older's inequality)}  \notag \\
&\leq &W^{p}\left( \int_{-\infty }^{\infty }\left( \int_{0}^{\infty }(\tilde{%
\gamma}(y,x))^{p}\frac{1}{\left( y^{2\alpha }+\xi ^{2\alpha }\right)
^{p\beta /2}}dy\right) \left( \int_{0}^{\infty }\frac{1}{\left( y^{2\alpha
}+\xi ^{2\alpha }\right) ^{q\beta /2}}dy\right) ^{p/q}dx\right)   \notag \\
&=&W^{p}\left( \frac{\Gamma \left( \frac{1}{2\alpha }\right) \Gamma \left( 
\frac{q\beta }{2}-\frac{1}{2\alpha }\right) }{2\Gamma \left( \frac{q\beta }{2%
}\right) \alpha \xi ^{q\alpha \beta -1}}\right) ^{p/q}\left( \int_{-\infty
}^{\infty }\left( \int_{0}^{\infty }(\tilde{\gamma}(y,x))^{p}\frac{1}{\left(
y^{2\alpha }+\xi ^{2\alpha }\right) ^{p\beta /2}}dy\right) dx\right)   \notag
\\
&=&\frac{\left[ \Gamma \left( \beta \right) \right] ^{p}\alpha \xi ^{\alpha
\beta p-1}\Gamma \left( \frac{q\beta }{2}-\frac{1}{2\alpha }\right) ^{p/q}}{%
2^{\frac{p}{q}}\Gamma \left( \frac{1}{2\alpha }\right) \Gamma \left( \beta -%
\frac{1}{2\alpha }\right) ^{p}\Gamma \left( \frac{q\beta }{2}\right) ^{p/q}}%
\left( \int_{-\infty }^{\infty }\left( \int_{0}^{\infty }(\tilde{\gamma}%
(y,x))^{p}\frac{1}{\left( y^{2\alpha }+\xi ^{2\alpha }\right) ^{p\beta /2}}%
dy\right) dx\right)   \notag \\
&=&:(\ast ).
\end{eqnarray}

By applying again H\"{o}lder's inequality we see that 
\begin{equation}
\tilde{\gamma}(y,x)\leq \frac{\left( \int_{0}^{y}|\tilde{\Delta}%
_{t}^{2}f^{(n)}(x)|^{p}dt\right) ^{1/p}}{(n-1)!}\frac{y^{(n-1+\frac{1}{q})}}{%
(q(n-1)+1)^{1/q}}\,.
\end{equation}%
Therefore it holds 
\begin{eqnarray}
(\ast ) &\leq &\frac{\left[ \Gamma \left( \beta \right) \right] ^{p}\alpha
\xi ^{\alpha \beta p-1}\Gamma \left( \frac{q\beta }{2}-\frac{1}{2\alpha }%
\right) ^{p/q}}{2^{\frac{p}{q}}\Gamma \left( \frac{1}{2\alpha }\right)
\Gamma \left( \beta -\frac{1}{2\alpha }\right) ^{p}\Gamma \left( \frac{%
q\beta }{2}\right) ^{p/q}\left[ (n-1)!\right] ^{p}(q(n-1)+1)^{p/q}}  \notag
\\
&&\cdot \left( \int_{-\infty }^{\infty }\left( \int_{0}^{\infty }\left(
\int_{0}^{y}|\tilde{\Delta}_{t}^{2}f^{(n)}(x)|^{p}dt\right) y^{p(n-1+\frac{1%
}{q})}\frac{1}{\left( y^{2\alpha }+\xi ^{2\alpha }\right) ^{p\beta /2}}%
dy\right) dx\right)  \notag \\
&=&\frac{\left[ \Gamma \left( \beta \right) \right] ^{p}\alpha \xi ^{\alpha
\beta p-1}\Gamma \left( \frac{q\beta }{2}-\frac{1}{2\alpha }\right) ^{p/q}}{%
2^{\frac{p}{q}}\Gamma \left( \frac{1}{2\alpha }\right) \Gamma \left( \beta -%
\frac{1}{2\alpha }\right) ^{p}\Gamma \left( \frac{q\beta }{2}\right) ^{p/q}%
\left[ (n-1)!\right] ^{p}(q(n-1)+1)^{p/q}}  \notag \\
&&\cdot \left( \int_{0}^{\infty }\left( \int_{-\infty }^{\infty }\left(
\int_{0}^{y}|\tilde{\Delta}_{t}^{2}f^{(n)}(x)|^{p}dt\right) y^{p(n-1+\frac{1%
}{q})}\frac{1}{\left( y^{2\alpha }+\xi ^{2\alpha }\right) ^{p\beta /2}}%
dx\right) dy\right)  \notag \\
&=&:(\ast \ast ).
\end{eqnarray}%
We call 
\begin{equation}
c_{2}:=\frac{\left[ \Gamma \left( \beta \right) \right] ^{p}\alpha \xi
^{\alpha \beta p-1}\Gamma \left( \frac{q\beta }{2}-\frac{1}{2\alpha }\right)
^{p/q}}{2^{\frac{p}{q}}\Gamma \left( \frac{1}{2\alpha }\right) \Gamma \left(
\beta -\frac{1}{2\alpha }\right) ^{p}\Gamma \left( \frac{q\beta }{2}\right)
^{p/q}\left[ (n-1)!\right] ^{p}(q(n-1)+1)^{p/q}}\,.
\end{equation}%
And hence 
\begin{eqnarray}
(\ast \ast ) &=&c_{2}\left( \int_{0}^{\infty }\left( \int_{-\infty }^{\infty
}\left( \int_{0}^{y}|\tilde{\Delta}_{t}^{2}f^{(n)}(x)|^{p}dt\right)
dx\right) y^{p(n-1+\frac{1}{q})}\frac{1}{\left( y^{2\alpha }+\xi ^{2\alpha
}\right) ^{p\beta /2}}dy\right)  \notag \\
&=&c_{2}\left( \int_{0}^{\infty }\left( \int_{0}^{y}\left( \int_{-\infty
}^{\infty }|\tilde{\Delta}_{t}^{2}f^{(n)}(x)|^{p}dx\right) dt\right) y^{pn-1}%
\frac{1}{\left( y^{2\alpha }+\xi ^{2\alpha }\right) ^{p\beta /2}}dy\right) 
\notag \\
&=&c_{2}\left( \int_{0}^{\infty }\left( \int_{0}^{y}\left( \int_{-\infty
}^{\infty }|\Delta _{t}^{2}f^{(n)}(x-t)|^{p}dx\right) dt\right) y^{pn-1}%
\frac{1}{\left( y^{2\alpha }+\xi ^{2\alpha }\right) ^{p\beta /2}}dy\right) 
\notag \\
&=&c_{2}\left( \int_{0}^{\infty }\left( \int_{0}^{y}\left( \int_{-\infty
}^{\infty }|\Delta _{t}^{2}f^{(n)}(x)|^{p}dx\right) dt\right) y^{pn-1}\frac{1%
}{\left( y^{2\alpha }+\xi ^{2\alpha }\right) ^{p\beta /2}}dy\right)  \notag
\\
&\leq &c_{2}\left( \int_{0}^{\infty }\left( \int_{0}^{y}\omega
_{2}(f^{(n)},t)_{p}^{p}dt\right) y^{pn-1}\frac{1}{\left( y^{2\alpha }+\xi
^{2\alpha }\right) ^{p\beta /2}}dy\right)  \notag \\
&\leq &c_{2}\omega _{2}(f^{(n)},\xi )_{p}^{p}\left( \int_{0}^{\infty }\left(
\int_{0}^{y}\left( 1+\frac{t}{\xi }\right) ^{2p}dt\right) y^{pn-1}\frac{1}{%
\left( y^{2\alpha }+\xi ^{2\alpha }\right) ^{p\beta /2}}dy\right) .
\end{eqnarray}%
I.e. so far we proved that 
\begin{equation}
\Lambda \leq c_{2}\omega _{2}(f^{(n)},\xi )_{p}^{p}\left( \int_{0}^{\infty
}\left( \int_{0}^{y}\left( 1+\frac{t}{\xi }\right) ^{2p}dt\right) y^{pn-1}%
\frac{1}{\left( y^{2\alpha }+\xi ^{2\alpha }\right) ^{p\beta /2}}dy\right) .
\end{equation}%
But 
\begin{equation}
\hbox{R.H.S.}(68)=\frac{c_{2}\xi }{2p+1}\omega _{2}(f^{(n)},\xi
)_{p}^{p}\left( \int_{0}^{\infty }\left( \left( 1+\frac{y}{\xi }\right)
^{2p+1}-1\right) y^{pn-1}\frac{1}{\left( y^{2\alpha }+\xi ^{2\alpha }\right)
^{p\beta /2}}dy\right) .
\end{equation}%
Call 
\begin{equation}
M:=\int_{0}^{\infty }\left( \left( 1+\frac{y}{\xi }\right) ^{2p+1}-1\right)
y^{pn-1}\frac{1}{\left( y^{2\alpha }+\xi ^{2\alpha }\right) ^{p\beta /2}}dy,
\end{equation}%
and 
\begin{equation}
\tilde{\tau}:=\int_{0}^{\infty }\left( \left( 1+u\right) ^{2p+1}-1\right)
u^{pn-1}\frac{1}{\left( 1+u^{2\alpha }\right) ^{p\beta /2}}du.
\end{equation}%
That is 
\begin{equation}
M=\xi ^{p\left( n-\alpha \beta \right) }\tilde{\tau}.
\end{equation}%
Therefore it holds 
\begin{equation}
\Lambda \leq \frac{\tilde{\tau}\left[ \Gamma \left( \beta \right) \right]
^{p}\alpha \xi ^{pn}\Gamma \left( \frac{q\beta }{2}-\frac{1}{2\alpha }%
\right) ^{p/q}\omega _{2}(f^{(n)},\xi )_{p}^{p}}{2^{\frac{p}{q}}\left(
2p+1\right) \Gamma \left( \frac{1}{2\alpha }\right) \Gamma \left( \beta -%
\frac{1}{2\alpha }\right) ^{p}\Gamma \left( \frac{q\beta }{2}\right) ^{p/q}%
\left[ (n-1)!\right] ^{p}(q(n-1)+1)^{p/q}}\,.
\end{equation}%
We have established (58). \hfill $\blacksquare $

\medskip The counterpart of Theorem 3 follows, $p =1$ case.

\medskip \noindent \textbf{Theorem 4.} \textit{Let } $f\in C^{n}(\mathbb{R})$
and $f^{(n)}\in L_{1}(\mathbb{R}),$\textit{\ }$n\geq 2$\textit{\ even, }$%
\alpha \in \mathbb{N},$ $\beta >\frac{n+3}{2\alpha }$\textit{. Then} 
\begin{eqnarray}
\Vert K(x)\Vert _{1} &\leq &\frac{1}{6\Gamma \left( \frac{1}{2\alpha }%
\right) \Gamma \left( \beta -\frac{1}{2\alpha }\right) (n-1)!}\left[ 3\Gamma
\left( \frac{n+1}{2\alpha }\right) \Gamma \left( \beta -\frac{n+1}{2\alpha }%
\right) \right. \\
&&\left. +3\Gamma \left( \frac{n+2}{2\alpha }\right) \Gamma \left( \beta -%
\frac{n+2}{2\alpha }\right) +\Gamma \left( \frac{n+3}{2\alpha }\right)
\Gamma \left( \beta -\frac{n+3}{2\alpha }\right) \right] \omega
_{2}(f^{(n)},\xi )_{1}\xi ^{n}.  \notag
\end{eqnarray}%
\textit{Hence as} $\xi \rightarrow 0$\textit{\ we obtain} $\Vert K(x)\Vert
_{1}\rightarrow 0$.

If additionally $f^{(2m)}\in L_{1}(\mathbb{R}),m=1,2,\ldots ,\frac{n}{2}$
then $\left\Vert M_{\xi }(f)-f\right\Vert _{1}\rightarrow 0,$ as $\xi
\rightarrow 0.$

\medskip \noindent \textbf{Proof.} Notice that 
\begin{equation}
\tilde{\Delta}^2_t f^{(n)} (x) =\Delta^2_t f^{(n)} (x -t),
\end{equation}
all $x,t\in \mathbb{R}$. Also it holds 
\begin{equation}
\int^\infty_{-\infty} |\Delta^2_t f^{(n)} (x -t)| dx = \int^\infty_{-\infty}
|\Delta^2_t f^{(n)} (w)| dw \leq \omega_2 (f^{(n)} ,t)_1 ,\quad \hbox{all} \
t\in \mathbb{R}_+ .
\end{equation}

Here we obtain%
\begin{eqnarray*}
\Vert K(x)\Vert _{1} &=&\int_{-\infty }^{\infty }|K(x)|dx \\
&&\overset{(57)}{\leq }W\int_{-\infty }^{\infty }\left( \int_{0}^{\infty
}\left( \int_{0}^{y}\left\vert \tilde{\Delta}_{t}^{2}f^{(n)}(x)\right\vert 
\frac{(y-t)^{n-1}}{(n-1)!}dt\right) \frac{1}{\left( y^{2\alpha }+\xi
^{2\alpha }\right) ^{\beta }}dy\right) dx \\
&\leq &W\int_{-\infty }^{\infty }\left( \int_{0}^{\infty }\left( \frac{%
y^{n-1}}{(n-1)!}\int_{0}^{y}\left\vert \tilde{\Delta}_{t}^{2}f^{(n)}(x)%
\right\vert dt\right) \frac{1}{\left( y^{2\alpha }+\xi ^{2\alpha }\right)
^{\beta }}dy\right) dx \\
&=&W\int_{0}^{\infty }\left( \left( \int_{-\infty }^{\infty }\left(
\int_{0}^{y}\left\vert \tilde{\Delta}_{t}^{2}f^{(n)}(x)\right\vert dt\right)
dx\right) \frac{y^{n-1}}{(n-1)!}\frac{1}{\left( y^{2\alpha }+\xi ^{2\alpha
}\right) ^{\beta }}\right) dy \\
&\overset{(75)}{=}&W\int_{0}^{\infty }\left( \left( \int_{-\infty }^{\infty
}\left( \int_{0}^{y}\left\vert \Delta _{t}^{2}f^{(n)}(x-t)\right\vert
dt\right) dx\right) \frac{y^{n-1}}{(n-1)!}\frac{1}{\left( y^{2\alpha }+\xi
^{2\alpha }\right) ^{\beta }}\right) dy \\
&=&W\int_{0}^{\infty }\left( \left( \int_{0}^{y}\left( \int_{-\infty
}^{\infty }\left\vert \Delta _{t}^{2}f^{(n)}(x-t)\right\vert dx\right)
dt\right) \frac{y^{n-1}}{(n-1)!}\frac{1}{\left( y^{2\alpha }+\xi ^{2\alpha
}\right) ^{\beta }}\right) dy \\
&&\overset{(76)}{\leq }W\int_{0}^{\infty }\left( \left( \int_{0}^{y}\omega
_{2}(f^{(n)},t)_{1}dt\right) \frac{y^{n-1}}{(n-1)!}\frac{1}{\left(
y^{2\alpha }+\xi ^{2\alpha }\right) ^{\beta }}\right) dy \\
&\leq &W\omega _{2}(f^{(n)},\xi )_{1}\left( \int_{0}^{\infty }\left(
\int_{0}^{y}\left( 1+\frac{t}{\xi }\right) ^{2}dt\right) \frac{y^{n-1}}{%
(n-1)!}\frac{1}{\left( y^{2\alpha }+\xi ^{2\alpha }\right) ^{\beta }}%
dy\right) \\
&=&W\omega _{2}(f^{(n)},\xi )_{1}\left( \int_{0}^{\infty }\left( \left( 1+%
\frac{y}{\xi }\right) ^{3}-1\right) \frac{\xi }{3}\frac{y^{n-1}}{(n-1)!}%
\frac{1}{\left( y^{2\alpha }+\xi ^{2\alpha }\right) ^{\beta }}dy\right) \\
&=&\frac{\Gamma \left( \beta \right) \alpha \xi ^{n}}{\Gamma \left( \frac{1}{%
2\alpha }\right) \Gamma \left( \beta -\frac{1}{2\alpha }\right) (n-1)!3}%
\omega _{2}(f^{(n)},\xi )_{1}\left( \int_{0}^{\infty }\left( \left(
1+Y\right) ^{3}-1\right) Y^{n-1}\frac{1}{\left( Y^{2\alpha }+1\right)
^{\beta }}dY\right)
\end{eqnarray*}%
\begin{eqnarray}
&=&\frac{\Gamma \left( \beta \right) \alpha \xi ^{n}}{\Gamma \left( \frac{1}{%
2\alpha }\right) \Gamma \left( \beta -\frac{1}{2\alpha }\right) (n-1)!3}%
\omega _{2}(f^{(n)},\xi )_{1}\left( \int_{0}^{\infty }\left(
3Y+3Y^{2}+Y^{3}\right) Y^{n-1}\frac{1}{\left( Y^{2\alpha }+1\right) ^{\beta }%
}dY\right)  \notag \\
&=&\frac{\Gamma \left( \beta \right) \alpha \xi ^{n}}{\Gamma \left( \frac{1}{%
2\alpha }\right) \Gamma \left( \beta -\frac{1}{2\alpha }\right) (n-1)!3}%
\omega _{2}(f^{(n)},\xi )_{1}\left( \int_{0}^{\infty }\left(
3Y^{n}+3Y^{n+1}+Y^{n+2}\right) \frac{1}{\left( Y^{2\alpha }+1\right) ^{\beta
}}dY\right)  \notag \\
&=&\frac{1}{6\Gamma \left( \frac{1}{2\alpha }\right) \Gamma \left( \beta -%
\frac{1}{2\alpha }\right) (n-1)!}\left[ 3\Gamma \left( \frac{n+1}{2\alpha }%
\right) \Gamma \left( \beta -\frac{n+1}{2\alpha }\right) \right.  \notag \\
&&\left. +3\Gamma \left( \frac{n+2}{2\alpha }\right) \Gamma \left( \beta -%
\frac{n+2}{2\alpha }\right) +\Gamma \left( \frac{n+3}{2\alpha }\right)
\Gamma \left( \beta -\frac{n+3}{2\alpha }\right) \right] \omega
_{2}(f^{(n)},\xi )_{1}\xi ^{n}.
\end{eqnarray}%
We have proved (74). \hfill $\blacksquare $

\medskip The related case here of $n =0$ comes next.

\medskip \noindent \textbf{Proposition 3.} \textit{Let }$p,q>1$ \textit{such
that }$\frac{1}{p}+\frac{1}{q}=1,$ $\alpha \in \mathbb{N},$ $\beta >\frac{1}{%
\alpha }\left( 2+\frac{1}{p}\right) $ \textit{and the rest as above. Then} 
\begin{equation}
\Vert M_{\xi }(f)-f\Vert _{p}\leq \frac{\rho ^{1/p}\Gamma \left( \beta
\right) \alpha ^{1/p}\Gamma \left( \frac{q\beta }{2}-\frac{1}{2\alpha }%
\right) ^{1/q}}{2^{\frac{1}{q}}\Gamma \left( \frac{1}{2\alpha }\right)
^{1/p}\Gamma \left( \beta -\frac{1}{2\alpha }\right) \Gamma \left( \frac{%
q\beta }{2}\right) ^{1/q}}\omega _{2}(f,\xi )_{p},
\end{equation}%
\textit{where} 
\begin{equation}
0<\rho :=\int_{0}^{\infty }\left( 1+y\right) ^{2p}\frac{1}{\left( y^{2\alpha
}+1\right) ^{\beta p/2}}dy<\infty .
\end{equation}%
\textit{Hence as }$\xi \rightarrow 0$\textit{\ we obtain }$M_{\xi
}\rightarrow I$ \textit{in the }$L_{p}$\textit{\ norm, }$p>1$.

\medskip \noindent \textbf{Proof.} From (55) we get 
\begin{equation}
|M_{\xi }(f;x)-f(x)|^{p}\leq W^{p}\left( \int_{0}^{\infty }\left\vert \tilde{%
\Delta}_{y}^{2}f(x)\right\vert \frac{1}{\left( y^{2\alpha }+\xi ^{2\alpha
}\right) ^{\beta }}dy\right) ^{p}.
\end{equation}%
We then estimate%
\begin{eqnarray}
\int_{-\infty }^{\infty }|M_{\xi }(f;x)-f(x)|^{p}dx &\leq
&W^{p}\int_{-\infty }^{\infty }\left( \int_{0}^{\infty }\left\vert \tilde{%
\Delta}_{y}^{2}f(x)\right\vert \frac{1}{\left( y^{2\alpha }+\xi ^{2\alpha
}\right) ^{\beta }}dy\right) ^{p}dx  \notag \\
&=&W^{p}\int_{-\infty }^{\infty }\left( \int_{0}^{\infty }\left\vert \tilde{%
\Delta}_{y}^{2}f(x)\right\vert \frac{1}{\left( y^{2\alpha }+\xi ^{2\alpha
}\right) ^{\beta /2}}\frac{1}{\left( y^{2\alpha }+\xi ^{2\alpha }\right)
^{\beta /2}}dy\right) ^{p}dx  \notag
\end{eqnarray}%
\begin{eqnarray}
&\leq &W^{p}\int_{-\infty }^{\infty }\left( \left( \int_{0}^{\infty }\left(
\left\vert \tilde{\Delta}_{y}^{2}f(x)\right\vert \frac{1}{\left( y^{2\alpha
}+\xi ^{2\alpha }\right) ^{\beta /2}}\right) ^{p}dy\right) ^{\frac{1}{p}%
}\left( \int_{0}^{\infty }\left( \frac{1}{\left( y^{2\alpha }+\xi ^{2\alpha
}\right) ^{\beta /2}}\right) ^{q}dy\right) ^{\frac{1}{q}}\right) ^{p}dx 
\notag \\
&=&W^{p}\int_{-\infty }^{\infty }\left( \left( \int_{0}^{\infty }\left\vert 
\tilde{\Delta}_{y}^{2}f(x)\right\vert ^{p}\frac{1}{\left( y^{2\alpha }+\xi
^{2\alpha }\right) ^{\beta p/2}}dy\right) \left( \int_{0}^{\infty }\frac{1}{%
\left( y^{2\alpha }+\xi ^{2\alpha }\right) ^{\beta q/2}}dy\right) ^{\frac{p}{%
q}}\right) dx  \notag \\
&=&W^{p}\left( \int_{0}^{\infty }\left( \int_{-\infty }^{\infty }\left\vert 
\tilde{\Delta}_{y}^{2}f(x)\right\vert ^{p}\frac{1}{\left( y^{2\alpha }+\xi
^{2\alpha }\right) ^{\beta p/2}}dx\right) dy\right) \left( \frac{\Gamma
\left( \frac{1}{2\alpha }\right) \Gamma \left( \frac{q\beta }{2}-\frac{1}{%
2\alpha }\right) }{2\Gamma \left( \frac{q\beta }{2}\right) \alpha \xi
^{q\alpha \beta -1}}\right) ^{\frac{p}{q}}  \notag
\end{eqnarray}%
\begin{eqnarray}
&=&\frac{\left[ \Gamma \left( \beta \right) \right] ^{p}\alpha \xi ^{\alpha
\beta p-1}\Gamma \left( \frac{q\beta }{2}-\frac{1}{2\alpha }\right) ^{p/q}}{%
2^{\frac{p}{q}}\Gamma \left( \frac{1}{2\alpha }\right) \Gamma \left( \beta -%
\frac{1}{2\alpha }\right) ^{p}\Gamma \left( \frac{q\beta }{2}\right) ^{p/q}}
\notag \\
&&\cdot \left( \int_{0}^{\infty }\left( \int_{-\infty }^{\infty }\left\vert
\Delta _{y}^{2}f(x-y)\right\vert ^{p}dx\right) \frac{1}{\left( y^{2\alpha
}+\xi ^{2\alpha }\right) ^{\beta p/2}}dy\right)   \notag \\
&\leq &\frac{\left[ \Gamma \left( \beta \right) \right] ^{p}\alpha \xi
^{\alpha \beta p-1}\Gamma \left( \frac{q\beta }{2}-\frac{1}{2\alpha }\right)
^{p/q}}{2^{\frac{p}{q}}\Gamma \left( \frac{1}{2\alpha }\right) \Gamma \left(
\beta -\frac{1}{2\alpha }\right) ^{p}\Gamma \left( \frac{q\beta }{2}\right)
^{p/q}}\int_{0}^{\infty }\omega _{2}(f,y)_{p}^{p}\frac{1}{\left( y^{2\alpha
}+\xi ^{2\alpha }\right) ^{\beta p/2}}dy  \notag \\
&\leq &\frac{\left[ \Gamma \left( \beta \right) \right] ^{p}\alpha \xi
^{\alpha \beta p-1}\Gamma \left( \frac{q\beta }{2}-\frac{1}{2\alpha }\right)
^{p/q}}{2^{\frac{p}{q}}\Gamma \left( \frac{1}{2\alpha }\right) \Gamma \left(
\beta -\frac{1}{2\alpha }\right) ^{p}\Gamma \left( \frac{q\beta }{2}\right)
^{p/q}}\omega _{2}(f,\xi )_{p}^{p}\int_{0}^{\infty }\left( 1+\frac{y}{\xi }%
\right) ^{2p}\frac{1}{\left( y^{2\alpha }+\xi ^{2\alpha }\right) ^{\beta p/2}%
}dy  \notag \\
&=&\frac{\left[ \Gamma \left( \beta \right) \right] ^{p}\alpha \xi ^{\alpha
\beta p-1}\Gamma \left( \frac{q\beta }{2}-\frac{1}{2\alpha }\right) ^{p/q}}{%
2^{\frac{p}{q}}\Gamma \left( \frac{1}{2\alpha }\right) \Gamma \left( \beta -%
\frac{1}{2\alpha }\right) ^{p}\Gamma \left( \frac{q\beta }{2}\right) ^{p/q}}%
\omega _{2}(f,\xi )_{p}^{p}\int_{0}^{\infty }\left( 1+Y\right) ^{2p}\frac{1}{%
\left( Y^{2\alpha }+1\right) ^{\beta p/2}\xi ^{\alpha \beta p}}\xi dY  \notag
\\
&=&\frac{\left[ \Gamma \left( \beta \right) \right] ^{p}\alpha \Gamma \left( 
\frac{q\beta }{2}-\frac{1}{2\alpha }\right) ^{p/q}}{2^{\frac{p}{q}}\Gamma
\left( \frac{1}{2\alpha }\right) \Gamma \left( \beta -\frac{1}{2\alpha }%
\right) ^{p}\Gamma \left( \frac{q\beta }{2}\right) ^{p/q}}\omega _{2}(f,\xi
)_{p}^{p}\int_{0}^{\infty }\left( 1+y\right) ^{2p}\frac{1}{\left( y^{2\alpha
}+1\right) ^{\beta p/2}}dy.
\end{eqnarray}%
The proof of (78) is now completed. \hfill $\blacksquare $

\medskip Also we give

\medskip \noindent \textbf{Proposition 4.} \textit{\ For }$\alpha \in 
\mathbb{N},$ $\beta >\frac{3}{2\alpha },$\textit{\ it holds,} 
\begin{equation}
\Vert M_{\xi }f-f\Vert _{1}\leq \left[ \frac{1}{2}+\frac{\Gamma \left( \frac{%
1}{\alpha }\right) \Gamma \left( \beta -\frac{1}{\alpha }\right) }{\Gamma
\left( \frac{1}{2\alpha }\right) \Gamma \left( \beta -\frac{1}{2\alpha }%
\right) }+\frac{\Gamma \left( \frac{3}{2\alpha }\right) \Gamma \left( \beta -%
\frac{3}{2\alpha }\right) }{2\Gamma \left( \frac{1}{2\alpha }\right) \Gamma
\left( \beta -\frac{1}{2\alpha }\right) }\right] \omega _{2}(f,\xi )_{1}.
\end{equation}%
\textit{Hence as} $\xi \rightarrow 0$\textit{\ we get }$M_{\xi }\rightarrow
I $\textit{\ in the} $L_{1}$\textit{\ norm.}

\medskip \noindent \textbf{Proof.} From (55) we have 
\begin{equation}
|M_{\xi }(f;x)-f(x)|\leq W\left( \int_{0}^{\infty }\left\vert \tilde{\Delta}%
_{y}^{2}f(x)\right\vert \frac{1}{\left( y^{2\alpha }+\xi ^{2\alpha }\right)
^{\beta }}dy\right) .
\end{equation}%
Hence we get%
\begin{eqnarray}
\int_{-\infty }^{\infty }|M_{\xi }(f;x)-f(x)|dx &\leq &W\int_{-\infty
}^{\infty }\left( \int_{0}^{\infty }|\tilde{\Delta}_{y}^{2}f(x)|\frac{1}{%
\left( y^{2\alpha }+\xi ^{2\alpha }\right) ^{\beta }}dy\right) dx  \notag \\
&=&W\int_{0}^{\infty }\left( \int_{-\infty }^{\infty }|\tilde{\Delta}%
_{y}^{2}f(x)|dx\right) \frac{1}{\left( y^{2\alpha }+\xi ^{2\alpha }\right)
^{\beta }}dy  \notag \\
&=&W\int_{0}^{\infty }\left( \int_{-\infty }^{\infty }|\Delta
_{y}^{2}f(x-y)|dx\right) \frac{1}{\left( y^{2\alpha }+\xi ^{2\alpha }\right)
^{\beta }}dy  \notag \\
\qquad &=&\;W\int_{0}^{\infty }\left( \int_{-\infty }^{\infty }|\Delta
_{y}^{2}f(x)|dx\right) \frac{1}{\left( y^{2\alpha }+\xi ^{2\alpha }\right)
^{\beta }}dy  \notag \\
\qquad &\leq &\;W\int_{0}^{\infty }\omega _{2}(f,y)_{1}\frac{1}{\left(
y^{2\alpha }+\xi ^{2\alpha }\right) ^{\beta }}dy  \notag
\end{eqnarray}
\begin{eqnarray}
&\leq &W\omega _{2}(f,\xi )_{1}\int_{0}^{\infty }\left( 1+\frac{y}{\xi }%
\right) ^{2}\frac{1}{\left( y^{2\alpha }+\xi ^{2\alpha }\right) ^{\beta }}dy
\notag \\
&=&W\omega _{2}(f,\xi )_{1}\int_{0}^{\infty }\left( 1+x\right) ^{2}\frac{1}{%
\left( x^{2\alpha }+1\right) ^{\beta }\xi ^{2\alpha \beta }}\xi dx  \notag \\
&=&\frac{\Gamma \left( \beta \right) \alpha }{\Gamma \left( \frac{1}{2\alpha 
}\right) \Gamma \left( \beta -\frac{1}{2\alpha }\right) }\omega _{2}(f,\xi
)_{1}\int_{0}^{\infty }\left( 1+x\right) ^{2}\frac{1}{\left( x^{2\alpha
}+1\right) ^{\beta }}dx  \notag \\
&=&\left[ \frac{1}{2}+\frac{\Gamma \left( \frac{1}{\alpha }\right) \Gamma
\left( \beta -\frac{1}{\alpha }\right) }{\Gamma \left( \frac{1}{2\alpha }%
\right) \Gamma \left( \beta -\frac{1}{2\alpha }\right) }+\frac{\Gamma \left( 
\frac{3}{2\alpha }\right) \Gamma \left( \beta -\frac{3}{2\alpha }\right) }{%
2\Gamma \left( \frac{1}{2\alpha }\right) \Gamma \left( \beta -\frac{1}{%
2\alpha }\right) }\right] \omega _{2}(f,\xi )_{1}.
\end{eqnarray}%
We have established (82).\hfill $\blacksquare $

\section*{References}

\begin{enumerate}
\item[{[1]}] G.A. Anastassiou, Rate of convergence of non-positive linear
convolution type operators. A sharp inequality, \textit{J. Math. Anal. and
Appl.}, \textbf{142} (1989), 441--451.

\item[{[2]}] G.A. Anastassiou, Sharp inequalities for convolution type
operators, \textit{Journal of Approximation Theory}, \textbf{58} (1989),
259--266.

\item[{[3]}] G.A. Anastassiou, \textit{Moments in Probability and
Approximation Theory}, Pitman Research Notes in Math., Vol. 287, Longman
Sci. \& Tech., Harlow, U.K., 1993.

\item[{[4]}] G.A. Anastassiou, \textit{Quantitative Approximations}, Chapman
\& Hall/CRC, Boca Raton, New York, 2001.

\item[{[5]}] G.A. Anastassiou, Basic convergence with rates of smooth Picard
singular integral operators, \textit{J. of Computational Analysis and
Applications},Vol.8, No.4 (2006), 313-334.

\item[{[6]}] G.A. Anastassiou, L$_{p}$ convergence with rates of smooth
Picard singular operators, \textit{Differential \& difference equations and
applications, Hindawi Publ. Corp., New York, }(2006), 31--45.

\item[{[7]}] G.A. Anastassiou and S. Gal, Convergence of generalized singular
integrals to the unit, univariate case, \textit{Math. Inequalities \&
Applications}, \textbf{3}, No. 4 (2000), 511--518.

\item[{[8]}] G.A. Anastassiou and S. Gal, Convergence of generalized singular
integrals to the unit, multivariate case, \textit{Applied Math. Rev.}, Vol.
1, World Sci. Publ. Co., Singapore, 2000, pp. 1--8.

\item[{[9]}] G.A. Anastassiou and R. Mezei, \textit{L}$_{\mathit{p}}$\textit{%
\ Convergence with Rates of Smooth Gauss-Weierstrass Singular Operators},
Nonlinear Studies\textit{,} accepted 2008.

\item[{[10]}] G.A. Anastassiou and R. A. Mezei, \textit{Global Smoothness and
Uniform Convergence of Smooth Poisson-Cauchy Type Singular Operators,}
submitted 2009.

\item[{[11]}] G.A. Anastassiou and R. A. Mezei, \textit{A Voronovskaya Type
Theorem for Poisson-Cauchy Type Singular Operators,} submitted 2009.

\item[{[12]}] R.A. DeVore and G.G. Lorentz, \textit{Constructive Approximation%
}, Springer-Verlag, Vol. 303, Berlin, New York, 1993.

\item[{[13]}] S.G. Gal, Remark on the degree of approximation of continuous
functions by singular integrals, \textit{Math. Nachr.}, \textbf{164} (1993),
197--199.

\item[{[14]}] S.G. Gal, Degree of approximation of continuous functions by
some singular integrals, \textit{Rev. Anal. Num\'{e}r, Th\'{e}or. Approx.},
(Cluj), Tome XXVII, No. 2 (1998), 251--261.

\item[{[15]}] R.N. Mohapatra and R.S. Rodriguez, On the rate of convergence
of singular integrals for H\"{o}lder continuous functions, \textit{Math.
Nachr.} \textbf{149} (1990), 117--124.

\item[{[16]}] D. Zwillinger, \textit{CRC Standard Mathematical Tables and
Formulae, 30th Edition,} Chapman \& Hall/CRC, Boca Raton, 1995.
\end{enumerate}

\end{document}